\documentclass[a4, 11pt] {amsart} 
\usepackage{mathrsfs}
\usepackage{amsmath}
\usepackage{amssymb}
\usepackage[all]{xy}
\usepackage[dvips]{graphicx, color}
\usepackage{amsmath}
\usepackage{amssymb}
\usepackage[all]{xy}

\usepackage{graphicx}

\begin{document}

\numberwithin{equation}{section}
\newtheorem{thm}[equation]{Theorem}
\newtheorem{pro}[equation]{Proposition}
\newtheorem{prob}[equation]{Problem}
\newtheorem{qu}[equation]{Question}
\newtheorem{cor}[equation]{Corollary}
\newtheorem{con}[equation]{Conjecture}
\newtheorem{lem}[equation]{Lemma}
\theoremstyle{definition}
\newtheorem{ex}[equation]{Example}
\newtheorem{defn}[equation]{Definition}
\newtheorem{ob}[equation]{Observation}
\newtheorem{rem}[equation]{Remark}

\hyphenation{homeo-morphism} 

\newcommand{\calA}{\mathcal{A}} 
\newcommand{\calD}{\mathcal{D}} 
\newcommand{\calE}{\mathcal{E}}
\newcommand{\calC}{\mathcal{C}} 
\newcommand{\Set}{\mathcal{S}et\,} 
\newcommand{\Top}{\mathcal{T}\!op \,}
\newcommand{\Topst}{\mathcal{T}\!op\, ^*}
\newcommand{\calK}{\mathcal{K}} 
\newcommand{\calO}{\mathcal{O}} 
\newcommand{\calS}{\mathcal{S}} 
\newcommand{\calT}{\mathcal{T}} 
\newcommand{\Z}{{\mathbb Z}}
\newcommand{\C}{{\mathbb C}}
\newcommand{\Q}{{\mathbb Q}}
\newcommand{\R}{{\mathbb R}}
\newcommand{\N}{{\mathbb N}}
\newcommand{\F}{{\mathcal F}} 

\hfill

\title{Local comparisons of homological and homotopical mixed Hodge polynomials}  

\author{Shoji Yokura}
\thanks{
\\
}

\date{}

\address{Department of Mathematics  and Computer Science, Graduate School of Science and Engineering, Kagoshima University, 1-21-35 Korimoto, Kagoshima, 890-0065, Japan
}

\email{yokura@sci.kagoshima-u.ac.jp, k8607624@kadai.jp}

\maketitle

\begin{abstract} For a simply connected complex algebraic variey $X$, by the mixed Hodge structures $(W_{\bullet}, F^{\bullet})$ and $(\tilde W_{\bullet}, \tilde F^{\bullet})$ of the homology group $H_{*}(X;\mathbb Q)$ and the homotopy groups $\pi_{*}(X)\otimes \mathbb Q$ respectively, we have the following mixed Hodge polynomials
$$MH_X(t,u,v):= \sum_{k,p,q} \dim \Bigl ( Gr_{F^{\bullet}}^{p} Gr^{W_{\bullet}}_{p+q} H_k (X;\mathbb C)  \Bigr) t^{k} u^{-p} v^{-q},$$
$$\quad \, \, MH^{\pi}_X(t,u,v):= \sum_{k,p,q} \dim  \Bigl (Gr_{\tilde F^{\bullet}}^{p} Gr^{\tilde W_{\bullet}}_{p+q} (\pi_k(X) \otimes \mathbb C) \Bigr ) t^ku^{-p} v^{-q},$$
which are respectively called \emph{the homological mixed Hodge polynomial} and \emph{the homotopical mixed Hodge polynomial}. In this paper we discuss some inequalities concerning these two mixed Hodge polynomials. 
\end{abstract}

\section{Introduction}

For a complex algebraic variety $X$ there exists a mixed Hodge structure $(W_{\bullet}, F^{\bullet})$ on the homology group $H_{*}(X;\mathbb Q)$(\cite{De1, De2}). In \cite{Mo} J. W. Morgan first put
mixed Hodge structures on the rational homotopy groups in the smooth case. Then, Morgan's results were extended to singular varieties by R. M. Hain \cite{Hain 2} (cf. \cite{Hain 1}) and V. Navarro-Aznar \cite{Nav} independently (e.g., see 
\cite[p.234, Historical Remarks]{PS}).

Then as defined in the abstract we can define the following polynomials of three variables $t,u,v$ (see Remark \ref{rem1} below) :
$$MH_X(t,u,v) := \sum_{k,p,q} \dim \Bigl ( Gr_{F^{\bullet}}^{p} Gr^{W_{\bullet}}_{p+q} H_k (X;\mathbb C)  \Bigr) t^{k} u^{-p} v^{-q},$$
$$\quad \, \, MH^{\pi}_X(t,u,v) := \sum_{k,p,q} \dim  \Bigl (Gr_{\tilde F^{\bullet}}^{p} Gr^{\tilde W_{\bullet}}_{p+q} (\pi_k(X) \otimes \mathbb C) \Bigr ) t^ku^{-p} v^{-q}.$$

\begin{rem}\label{rem1} In this paper we consider the rational homology groups  $H_k(X;\mathbb Q)$ instead of the cohomology groups $H^k(X;\mathbb Q) \cong Hom(H_k(X;\mathbb Q), \mathbb Q)$ (by the universal coefficient theorem), thus the mixed Hodge structures have both $p, q$ negative, thus negative weights. Therefore in defining the mixed Hodge polynomial $MH_X(t,u,v)$ we consider $u^{-p} v^{-q}$ instead of $u^pv^q$ (cf. \cite[p.35]{PS}). It is the same for the homotopical mixed Hodge polynomial $MH^{\pi}_X(t,u,v)$. In other words, the above two polynomials can be also defined respectively using the cohomology groups $H^k (X;\mathbb C)$ and the dual $(\pi_k(X)\otimes \mathbb C)^{\vee}= Hom(\pi_k(X)\otimes \mathbb C; \mathbb C)$ of the homotopy group $\pi_k(X)\otimes \mathbb C$ by
$$MH_X(t,u,v) := \sum_{k,p,q} \dim \Bigl ( Gr_{F^{\bullet}}^{p} Gr^{W_{\bullet}}_{p+q} H^k (X;\mathbb C)  \Bigr) t^{k} u^p v^q, $$
$$\quad \quad \, \, \, MH^{\pi}_X(t,u,v) :=
\sum_{k,p,q} \dim  \Bigl (Gr_{\tilde F^{\bullet}}^{p} Gr^{\tilde W_{\bullet}}_{p+q} ((\pi_k(X)\otimes \mathbb C)^{\vee})\Bigr ) t^ku^p v^q.
$$
\end{rem}

\begin{rem} In order to get the mixed Hodge structure on the homotopy groups, in fact it suffices that the algebraic variety is nilpotent in the sense that $\pi_1$ is nilpotent and acting nilpotently on higher homotopy groups (e.g.,see \cite[Remark 8.12]{PS}). Simply connected is then a particular case.  
\end{rem}

The first polynomial is well-known, usually called the mixed Hodge polynomial and has been studied very well. The second one is a homotopical analogue, defined by the mixed Hodge structure on the homotopy groups $\pi_{*}(X)$. So, we call these two polynomials respectively \emph{the homological mixed Hodge polynomial} and \emph{the homotopical mixed Hodge polynomial}.

Here we observe the following for the special values $(u,v)=(1,1)$:\\

$$P_{X}(t) = MH_X(t,1,1) = \sum_{k\geqq 0} \dim H_k(X;\mathbb C) t^{k}  = 1 + \sum_{k\geqq 1} \dim H_k(X;\mathbb C) t^{k},$$
$$\, \, \, P^{\pi}_{X}(t) = MH^{\pi}_X(t,1,1) = \sum_{k\geqq 2} \dim (\pi_{k}(X) \otimes \mathbb C) t^{k} =\sum_{k\geqq 2} \dim (\pi_{k}(X) \otimes \mathbb Q) t^{k}.$$

The first polynomial is the usual \emph{Poincar\'e polynomial} and the second one is its homotopical analogue, called the \emph{homotopical Poincar\'e polynomial}.

In this note we discuss some inequalities concerning these two mixed Hodge polynomials $MH_X(t,u,v)$ and $MH^{\pi}_X(t,u,v)$. More details will appear elsewhere.

\section{Homological mixed Hodge polynomial and homotopical mixed Hodge polynomial}

The most important and fundamental topological invariant in geometry and topology is the Euler--Poincar\'e characteristic $\chi(X)$,
which is defined to be the alternating sum of the Betti numbers $\beta_i(X):=\dim_{\mathbb Q} H_i(X;\Q) = \dim_{\mathbb C} H_i(X;\mathbb C) $:
$$\chi(X):= \sum_{i \geqq 0} (-1)^i\beta_i(X),$$
provided that each $\beta_i(X)$ and $\chi(X)$ are both finite.
Similarly, for a topological space whose fundamental group is an Abelian group one can define the \emph{homotopical Betti number} $\beta^{\pi}_i(X):= \dim (\pi_i(X)\otimes \Q)$ where $i\geqq 1$ and the \emph{homotopical Euler--Poincar\'e characteristic}:
$$\chi^{\pi}(X):= \sum_{i \geqq 1}  (-1)^i\beta^{\pi}_i(X),$$
provided that each $\beta^{\pi}_i(X)$ and $\chi^{\pi}(X)$ are both finite.
The Euler--Poincar\'e characteristic is the special value of the Poincar\'e polynomial $P_X(t)$ at $t=-1$ and the homotopical Euler--Poincar\'e characteristic is the special value of the homotopical Poincar\'e polynomial $ P^{\pi}_X(t)$ at $t=-1$:
$$P_X(t):= \sum_{i \geqq 0}  t^i \beta_i(X), \quad \chi(X) = P_X(-1),$$
$$ P^{\pi}_X(t):= \sum_{i \geqq 1} t^i \beta^{\pi}_i(X), \quad \chi^{\pi}(X) = P^{\pi}_X(-1).$$

The Poincar\'e polynomial $P_X(t)$ is \emph{multiplicative} in the following sense:
$$P_{X \times Y}(t) = P_X(t) \times P_Y(t),$$
which follows from the K\"unneth Formula:
$$H_n(X \times Y;\Q)= \sum_{i+j=n} H_i(X; \Q) \otimes H_j(Y;\Q).$$
The homotopical Poincar\'e polynomial $P^{\pi}_X(t)$ is \emph{additive} in the following sense:
$$P^{\pi}_{X \times Y}(t) = P^{\pi}_X(t) + P^{\pi}_Y(t),$$
which follows from 
$$\pi_i(X \times Y) = \pi_i (X) \times \pi_i(Y) =\pi_i(X) \oplus \pi_i(Y)$$
and
$ (A \oplus B)\otimes \Q = (A \otimes \Q ) \oplus (B \otimes \Q).$

Here we note that
\begin{equation*}
P_X(t) = MH_X(t,1,1), \quad P^{\pi}_X(t) = MH^{\pi}_X(t,1,1).
\end{equation*}

In fact the homological mixed Hodge polynomial is also multiplicative just like the Poincar\'e polynomial $P_X(t)$
\begin{equation}\label{mh-multi}
MH_{X \times Y}(t,u,v) =MH_X(t,u,v)  \times MH_Y(t,u,v)
\end{equation}
which follows from the fact that the mixed Hodge structure is compatible with the tensor product (e.g., see \cite[\S 3.1, Examples 3.2]{PS}.)
As to the homotopical mixed Hodge polynomial, it is additive just like the homotopical Poincar\'e polynomial $P^{\pi}_X(t)$

\begin{equation}\label{mh-pi-additive}
MH^{\pi}_{X \times Y}(t,u,v) = MH^{\pi}_X(t,u,v) + MH^{\pi}_Y(t,u,v)
\end{equation}
since $\pi_{*}(X \times Y) = \pi_{*}(X) \oplus \pi_{*}(Y)$ and the category of mixed Hodge structures is abelian and the direct sum of a mixed Hodge structure is also a mixed Hodge structure. 
\section{Local comparisons of these two mixed Hodge polynomials}
By the above definition, we have $0= P^{\pi}_X(0) = MH^{\pi}_X(0,1,1) < MH_X(0,1,1) =P_X(0) = 1$. Hence we get the following strict inequality\footnote{We note that given two real valued polynomial (therefore, continuous) functions $f(x,y,z)$ and $g(x,y,z)$, a strict inequality $f(a,b,c) < g(a,b,c)$ at a special value $(a,b,c)$ implies a local strict inequality $f(x,y,z)<g(x,y,z)$ for $|x-a| \ll 1, |y-b| \ll 1, |z-c| \ll 1$}:

\begin{cor}\label{cor011}
$$MH^{\pi}_X(t,u,v) < MH_X(t,u,v)$$
for $|t| \ll 1, |u-1| \ll 1, |v-1| \ll 1$.
\end{cor}

When $t=-1$, $MH_X(-1,1,1)=P_X(-1) = \chi(X)$ is the Euler--Poincar\'e characteristic and $MH^{\pi}_X(-1,1,1)=P^{\pi}_X(-1) = \chi^{\pi}(X)$ is the homotopical Euler--Poincar\'e characteristic.  In this case we do have the following theorem due to  F\'elix--Halperin--Thomas \cite[Proposition 32.16]{FHT}:
\begin{thm} We have $\chi^{\pi}(X) < \chi(X)$, namely $MH^{\pi}_X(-1,1,1) < MH_X(-1,1,1)$.
\end{thm}

Hence we get the following strict inequality:
\begin{cor}\label{cor-111}
$$MH^{\pi}_X(t,u,v) < MH_X(t,u,v)$$
for $|t+1|\ll1, |u-1| \ll 1, |v-1| \ll 1 $.
\end{cor}

As to the case when $(t,u,v)=(1,1,1)$, we have 
$$\quad \quad MH_X(1,1,1)=P_X(1) = \sum_{k\geqq 0} \dim H_k(X;\mathbb C)= 1 + \sum_{k\geqq 1} \dim H_k(X;\mathbb C),$$
$$MH^{\pi}_X(1,1,1) = P^{\pi}_X(1) = \sum_{k\geqq 2} \dim (\pi_{k}(X) \otimes \mathbb C). \hspace{3cm} $$
For these integers we do have the following Hilali conjecture \cite{Hil}, which has been solved affirmatively for many spaces such as smooth complex projective varieties and symplectic manifolds (e.g. see \cite{BFMM, HM, HM2}), but still open:
\begin{con}[Hilali conjecture]  
$$P^{\pi}_X(1) \leqq P_X(1),$$
i.e., $MH^{\pi}_X(1,1,1) \leqq MH_X(1,1,1).$
\end{con}
\begin{rem} The inequality $\leqq$ in the Hilali conjecture cannot be replaced by the strict inequality $<$.
It follows from the minimal model of the de Rham algebra of $\mathbb P ^n$ that we have (see \cite[Example 9.9]{PS})
$$\pi_k(\mathbb P^{n}) \otimes \Q 
=\begin{cases} 
0 & \, k \not =2,2n+1\\
\Q & \, k=2, 2n+1.
\end{cases}
$$
In particular, in the case when $n=1$, we have
$$MH^{\pi}_{\mathbb P^1}(t,u,v) = t^2uv + t^3u^2v^2, \quad MH_{\mathbb P^1}(t,u,v) = 1 + t^2uv.$$

\noindent
So we have that $MH^{\pi}_X(1,1,1) = MH_X(1,1,1) =2$, i.e. $P^{\pi}_X(1) = P_X(1)=2$.
We also remark that in the case of (non-strict) inequality $MH^{\pi}_X(1,1,1) \leqq MH_X(1,1,1)$, unlike Corollary \ref{cor011} and Corollary \ref{cor-111} we cannot expect the following local inequality$$MH^{\pi}_X(t,u,v) \leqq MH_X(t,u,v)$$
for $|t-1| \ll 1, |u-1| \ll 1, |v-1| \ll 1$. Indeed, clearly the following does not hold:
$$MH^{\pi}_{\mathbb P^1}(t,1,1) = t^2 + t^3 \leqq 1 + t^2 = MH_{\mathbb P^1}(t,1,1)$$
for $|t-1| \ll 1.$
\end{rem}

However, using the multiplicativity of the Poincar\'e polynomial $P_X(t)$ and the additivity of the homotopical Poincar\'e polynomial $P^{\pi}_X(t)$, we can get the following theorem, which kind of says that the Hilali conjecture holds ``modulo product'' \cite{Yo}:
\begin{thm} \label{hilali-product} There exists a positive integer $n_{0}$ such that for $\forall n \geqq n_0$ the following strict inequality holds:
$$P^{\pi}_{X^n}(1) < P_{X^n}(1).$$
\end{thm}
Hence, since $P^{\pi}_{X^n}(1) < P_{X^n}(1)$ means $MH^{\pi}_{X^n}(1,1,1) < MH_{X^n}(1,1,1)$,  we have that 
\begin{equation}\label{ineq-111}
MH^{\pi}_{X^n}(1,1,1) < MH_{X^n}(1,1,1) \, \, \text{for} \, \, \forall n \geqq n_0
\end{equation} In fact we can get the following strict inequality, which, should be noted, does not follow straightforwardly from the above strict inequality (\ref{ineq-111}) and requires a bit of work:

\begin{cor} There exists a positive integer $n_{0}$ such that for $\forall n \geqq n_0$
$$MH^{\pi}_{X^n}(t,u,v) <  MH_{X^n}(t,u,v)$$
for $|t-1| \ll 1, |u-1| \ll 1, |v-1| \ll 1$.
\end{cor}

In fact, in a similar way, using the multiplicativity of the mixed Hodge polynomial, i.e., $(\ref{mh-multi})$ and the additivity of the homotopical mixed Hodge polynomial, i.e., $(\ref{mh-pi-additive})$, we can show the following theorem. Let $\mathbb R_{>0}$ be the set of positive real numbers.
\begin{thm} Let $(s,a,b) \in (\mathbb R_{>0})^3$. Then there exists a positive integer $n_{(s,a,b)}$ such that for 
$\forall n \geqq n_{(s,a,b)}$
the following strict inequality holds
$$MH^{\pi}_{X^n}(t,u,v) <  MH_{X^n}(t,u,v).$$
for $|t-s| \ll 1, |u-a| \ll 1, |v-b| \ll 1$.
\end{thm}

The following theorem follows from the above theorem and the compactness of the following compact cube $\mathscr C_{\varepsilon,r}$.
\begin{thm} Let $\varepsilon, r$ be positive real numbers such that $0<\varepsilon \ll 1$ and $\varepsilon < r$ and $\mathscr C_{\varepsilon,r}:=[\varepsilon, r] \times [\varepsilon,r] \times [\varepsilon,r] \subset (\mathbb R_{> 0})^3$ be a cube.
Then there exists a positive integer $n_{\varepsilon,r}$ such that for 
$\forall n \geqq n_{\varepsilon,r}$
the following strict inequality holds
$$MH^{\pi}_{X^n}(t,u,v) <  MH_{X^n}(t,u,v)$$
for $\forall (t,u,v) \in \mathscr C_{\varepsilon,r}.$ 
\end{thm}

We would like to pose the following conjecture:
\begin{con} Let $\varepsilon$ be a positive real number such that $0<\varepsilon \ll 1$. There exist a positive integer $n_0$ such that for $\forall n \geqq n_0$ the following strict inequality holds
$$MH^{\pi}_{X^n}(t,u,v) <  MH_{X^n}(t,u,v)$$
for $\forall (t,u,v) \in [\varepsilon, \infty)^3 \subset (\mathbb R_{> 0})^3$.
\end{con}
In the case when $u=v=1$, i.e., in the case of $P^{\pi}_X(t)$ and $P_X(t)$, we do have the following ``half-global'' version of Theorem \ref{hilali-product}:
\begin{thm} Let $\varepsilon$ be a positive real number such that $0<\varepsilon \ll 1$.
There exists a positive integer $n_{0}$ such that for $\forall n \geqq n_0$ the following strict inequality holds:
$$P^{\pi}_{X^n}(t) < P_{X^n}(t) \quad (\forall t \in [\varepsilon, \infty)).$$\
\end{thm}

{\bf Acknowledgements:} The author would like to thank the anonymous referee for his/her very useful comments and suggestions. The author also would like to thank  Anatoly Libgober for his interest in this work and various comments and suggestions on an earlier version of the paper. Their comments and suggestions improved the paper. In fact, a joint paper with A. Libgober is in preparation and will contain detailed proofs of results of the present paper as
well as calculations and further information about mixed Hodge
polynomials of elliptic spaces. This work is supported by JSPS KAKENHI Grant Number JP19K03468.

\end{document}